\documentclass[11pt]{amsart}

\issueinfo{00}{0}{Month}{Year}
\copyrightinfo{2008}{University of Calgary}
\pagespan{1}{2}
\PII{ISSN 1715-0868}

\usepackage{graphicx}
\usepackage{epstopdf}
\usepackage{amsthm,amsmath,amssymb,amsxtra,amscd,mathtools,amsfonts,times,hyperref,mathrsfs}
\usepackage{verbatim}
\usepackage{multirow}
\usepackage{multicol}
\usepackage{url}
\usepackage{algorithmic}
\usepackage{enumerate}
\usepackage{array}
\usepackage[super]{nth}

\newtheorem{theorem}{Theorem}[section]

\newtheorem*{example*}{Example}

\newcommand\numberthis{\addtocounter{equation}{1}\tag{\theequation}}

\newtheoremstyle{myexample}{3pt}{3pt}{\rmfamily}{}{\itshape}{:}{ }{\thmname{#1}\thmnumber{ #2}\thmnote{ (#3)}}
\theoremstyle{myexample}

\newtheoremstyle{myremark}{3pt}{3pt}{\rmfamily}{}{\itshape}{:}{ }{\thmname{#1}}
\theoremstyle{myremark}
\newtheorem{remark}[theorem]{Remark}
\newtheorem*{observation*}{Observation}

\newtheoremstyle{conjecture}{3pt}{3pt}{\itshape}{}{\bfseries}{.}{ }{\thmname{#1}\thmnote{ (#3)}}
\theoremstyle{conjecture}

\newtheorem*{question*}{Question}

\newtheorem{theorem*}{Theorem}

\numberwithin{equation}{section}

\newcounter{algorithm}
\renewcommand{\thealgorithm}{\thesection.\arabic{algorithm}}

\begin{document}

\title{Pfaff's Method Revisited}

\author{Aritram Dhar}
\address{Department of Mathematics, University of Florida, Gainesville FL 32611, USA}
\email{aritramdhar@ufl.edu}

%

\date{(date1), and in revised form (date2).}

\thanks{}

\begin{abstract}
In 1797, Pfaff gave a simple proof of a ${}_3F_2$ hypergeometric series summation formula which was much later reproved by Andrews in 1996. In the same paper, Andrews also proved other well-known hypergeometric identities using Pfaff's method. In this paper, we prove a number of terminating $q$-hypergeometric series-product identities using Pfaff's method thereby providing a detailed account of its wide applicability.
\end{abstract}

\maketitle

\section{Introduction}\label{s1}
Let $\mathbb{N}$ be the set of natural numbers. For complex variables $a$ and $q$, define the conventional $q$-shifted factorial by\\
\begin{align*}
(a;q)_n := \prod\limits_{k = 0}^{n-1}(1-aq^k)    
\end{align*}\\
for any $n\in\mathbb{N}$ and $(a;q)_0 = 1$. For $|q| < 1$, define\\
\begin{align*}
(a;q)_{\infty} := \lim\limits_{n\rightarrow\infty}(a;q)_n.
\end{align*}\\
We can also define $(a;q)_n$ for all real numbers $n$ by\\
\begin{align*}
(a;q)_n := \dfrac{(a;q)_{\infty}}{(aq^n;q)_{\infty}}.    
\end{align*}\\
The $q$-shifted factorial for negative integers can be defined by\\
\begin{align}\label{eq11}
(a;q)_{-n} := \dfrac{1}{(aq^{-n};q)_n} = \dfrac{(-q/a)^n}{(q/a;q)_n}q^{\binom{n}{2}}.    
\end{align}\\
\par Following Gasper and Rahman's text \cite{Gas-Rah04}, the ${}_{r+1}F_r$ hypergeometric series and the ${}_{r+1}\phi_r$ unilateral basic hypergeometric series are defined, respectively, by\\
\begin{equation}
\begin{multlined}\label{eq12}
{}_{r+1}F_r \left[
	\setlength\arraycolsep{2pt}
	\begin{matrix}
		a_0,\ldots,a_r \\
		\multicolumn{2}{c}{
			\begin{matrix}
				b_1,\ldots,b_r 	
			\end{matrix}}
	\end{matrix} \hspace{2pt}
; z \right] := \sum_{n=0}^{\infty}\dfrac{(a_0,\ldots,a_r)_n}{(1,b_1,\ldots,b_r)_n}z^n,
\end{multlined}    
\end{equation}\\
where
\begin{align*}
(a_1,a_2,\ldots,a_k)_n := \prod\limits_{i=1}^{k}\prod\limits_{j=0}^{n-1}(a_i+j),    
\end{align*}\\
and\\
\begin{equation}
\begin{multlined}\label{eq13}
{}_{r+1}\phi_r \left[
	\setlength\arraycolsep{2pt}
	\begin{matrix}
		a_0,\ldots,a_r \\
		\multicolumn{2}{c}{
			\begin{matrix}
				b_1,\ldots,b_r 	
			\end{matrix}}
	\end{matrix} \hspace{2pt}
;q, z \right] := \sum_{n=0}^{\infty}\dfrac{(a_0,\ldots,a_r;q)_n}{(q,b_1,\ldots,b_r;q)_n}z^n,\quad \lvert z\rvert<1,
\end{multlined}    
\end{equation}\\
where
\begin{align*}
(a_1,a_2,\ldots,a_k;q)_n := \prod_{i=1}^{k}(a_i;q)_n.    
\end{align*}\\
The series \eqref{eq13} is called a \textit{terminating series} if one of its numerator parameters is of the form $q^{-m}$ where $m\in\mathbb{N}\cup\{0\}$. The series \eqref{eq13} is called \textit{balanced} or \textit{Saalsch\"utzian} if\\
\begin{align*}
b_1b_2\ldots b_r = qa_0a_1\ldots a_r.    
\end{align*}\\
The series \eqref{eq13} is said to be \textit{well-poised} if the parameters satisfy the relations\\
\begin{align*}
qa_0 = a_1b_1 = a_2b_2 = \ldots = a_rb_r,    
\end{align*}\\
\textit{very-well-poised} if, in addition,\\
\begin{align*}
a_1 = q\sqrt{a_0},\, a_2 = -q\sqrt{a_0}.\\    
\end{align*}
In 1797, Pfaff \cite{Pfaff97} gave ``the simplest proof'' (as Andrews writes in \cite{And96}) of the following ${}_3F_2$ hypergeometric identity,\\
\begin{equation*}
{}_3F_2 \left[
	\setlength\arraycolsep{2pt}
	\begin{matrix}
		-n,\quad\,\,\,\,\, a,\quad\,\,\,\,\, b \\
		\multicolumn{2}{c}{
			\begin{matrix}
			  \,\,\,\,\, c,\,\,\, 1+a+b-n-c 	
			\end{matrix}}
	\end{matrix} \hspace{8pt}
;1 \right] = \dfrac{(c-a,c-b)_n}{(c,c-a-b)_n}.   
\end{equation*}\\
\par Pfaff's result was ignored until Saalsch\"utz rediscovered it in 1890. In the second part of a series of three papers (which he refers to as ``the Pfaff Trilogy''), Andrews \cite{And96} then recreated the proof of Pfaff's ${}_3F_2$ sum and also proved some other hypergeometric identities using Pfaff's method in the same paper, for instance, the terminating Kummer theorem, Bailey's theorem, the theorems of Dougall and Lakin, to name a few.\\
\par Towards the end of his paper, Andrews \cite{And96} proved a $q$-analog of Pfaff's original theorem which was due to Jackson \cite{Jack10}. He then concludes with a brief description saying Pfaff's method is mostly effective for balanced and well-poised hypergeometric series (and $q$-series). In this paper, we give proofs of a plethora of well-known terminating $q$-hypergeometric series-product identities using Pfaff's method.\\
\par We now present the statements of the main terminating $q$-hypergeometric identites which we prove in this paper.
\\
\begin{theorem}($q$-binomial theorem \cite[II.$4$]{Gas-Rah04})\label{thm11}
For any non-negative integer $n$, we have\\
\begin{equation}\label{eq14}
{}_1\phi_0 \left[
	\setlength\arraycolsep{2pt}
	\begin{matrix}
		q^{-n} \\
		\multicolumn{2}{c}{
			\begin{matrix}
				- 	
			\end{matrix}}
	\end{matrix} \hspace{8pt}
;q, z \right] = (zq^{-n};q)_n.   
\end{equation}\\
\end{theorem}

\begin{theorem}($q$-Chu-Vandermonde \cite[II.$6$]{Gas-Rah04})\label{thm12}
For any non-negative integer $n$, we have\\
\begin{equation}\label{eq15}
{}_2\phi_1 \left[
	\setlength\arraycolsep{2pt}
	\begin{matrix}
		a, & q^{-n} \\
		\multicolumn{2}{c}{
			\begin{matrix}
				 c 	
			\end{matrix}}
	\end{matrix} \hspace{2pt}
;q, q \right] = \dfrac{(c/a;q)_n}{(c;q)_n}a^n.   
\end{equation}\\
\end{theorem}

\begin{theorem}($q$-Pfaff-Saalsch\"utz \cite[II.$12$]{Gas-Rah04})\label{thm13}
For any non-negative integer $n$, we have\\
\begin{equation}\label{eq16}
{}_3\phi_2 \left[
	\setlength\arraycolsep{2pt}
	\begin{matrix}
		a,\quad\,\,\,\,\, b,\quad\,\,\,\,\, q^{-n} \\
		\multicolumn{2}{c}{
			\begin{matrix}
			  \,\,\,\,\, c,\,\,\, abq^{-n+1}/c 	
			\end{matrix}}
	\end{matrix} \hspace{8pt}
;q, q \right] = \dfrac{(c/a,c/b;q)_n}{(c,c/ab;q)_n}.   
\end{equation}\\
\end{theorem}

\begin{theorem}($q$-Dixon \cite[II.$14$]{Gas-Rah04})\label{thm14}
For any non-negative integer $n$, we have\\
\begin{equation}\label{eq17}
{}_4\phi_3 \left[
	\setlength\arraycolsep{2pt}
	\begin{matrix}
		a,\quad -q\sqrt{a},\quad b,\quad q^{-n} \\
		\multicolumn{2}{c}{
			\begin{matrix}
			  \,\,\,\,\, -\sqrt{a},\,\,\, aq/b,\,\,\, aq^{n+1} 	
			\end{matrix}}
	\end{matrix} \hspace{8pt}
;q, \dfrac{q^{n+1}\sqrt{a}}{b} \right] = \dfrac{(aq,q\sqrt{a}/b;q)_n}{(q\sqrt{a},aq/b;q)_n}.   
\end{equation}\\
\end{theorem}

\begin{theorem}(\cite[p. $110$, Ex. $3.34$]{Gas-Rah04})\label{thm15}
For any non-negative integer $n$, we have\\
\begin{equation}\label{eq18}
{}_4\phi_3 \left[
	\setlength\arraycolsep{2pt}
	\begin{matrix}
		q^{-2n},\quad c^2,\quad a,\quad aq \\
		\multicolumn{2}{c}{
			\begin{matrix}
			  \,\,\,\,\, a^2q^2,\,\,\, cq^{-n},\,\,\, cq^{-n+1} 	
			\end{matrix}}
	\end{matrix} \hspace{8pt}
;q^2, q^2 \right] = \dfrac{(-q,aq/c;q)_n}{(-aq,q/c;q)_n}.   
\end{equation}\\
\end{theorem}

\begin{theorem}(Andrews \cite[p. $22$, eq.$(7.7)$ and eq.$(7.6)$]{And96}, Li-Chu \cite[p.$1006$, eq.$(2)$ and eq.$(3)$]{Li-Chu19})\label{thm16}
For any non-negative integer $n$, we have\\
\begin{equation}\label{eq19}
{}_4\phi_3 \left[
	\setlength\arraycolsep{2pt}
	\begin{matrix}
		q^{-n},\quad b,\quad b\sqrt{q},\quad d^2q^n \\
		\multicolumn{2}{c}{
			\begin{matrix}
			  \quad dq,\quad d\sqrt{q},\quad b^2 	
			\end{matrix}}
	\end{matrix} \hspace{28pt}
;q, q \right] = \dfrac{b^n(1-d)(-\sqrt{q},d\sqrt{q}/b;\sqrt{q})_n}{(1-q^nd)(-b,d;\sqrt{q})_n}   
\end{equation}\\
and\\
\begin{equation}\label{eq110}
{}_4\phi_3 \left[
	\setlength\arraycolsep{2pt}
	\begin{matrix}
		q^{-n},\quad b,\quad b\sqrt{q},\quad d^2q^{n+1} \\
		\multicolumn{2}{c}{
			\begin{matrix}
			  \quad dq,\quad d\sqrt{q},\quad b^2q 	
			\end{matrix}}
	\end{matrix} \hspace{32pt}
;q, q \right] = \dfrac{b^n(-\sqrt{q},d\sqrt{q}/b;\sqrt{q})_n}{(-b\sqrt{q},d\sqrt{q};\sqrt{q})_n}.   
\end{equation}\\
\end{theorem}

\begin{theorem}(\cite[II.$21$]{Gas-Rah04})\label{thm17}
For any non-negative integer $n$, we have\\
\begin{equation}\label{eq111}
{}_6\phi_5 \left[
	\setlength\arraycolsep{2pt}
	\begin{matrix}
		a,\quad q\sqrt{a},\quad -q\sqrt{a},\quad b,\quad c,\quad q^{-n} \\
		\multicolumn{2}{c}{
			\begin{matrix}
			  \,\,\,\,\, \sqrt{a},\,\,\, -\sqrt{a},\,\,\, aq/b,\,\,\, aq/c,\,\,\, aq^{n+1} 	
			\end{matrix}}
	\end{matrix} \hspace{8pt}
;q, \dfrac{aq^{n+1}}{bc} \right] = \dfrac{(aq,aq/bc;q)_n}{(aq/b,aq/c;q)_n}.   
\end{equation}\\
\end{theorem}

\begin{theorem}(\cite[II.$22$]{Gas-Rah04})\label{thm18}
For any non-negative integer $n$, we have\\
\begin{align*}
&{}_8\phi_7 \left[
	\setlength\arraycolsep{2pt}
	\begin{matrix}
		a,\quad q\sqrt{a},\quad -q\sqrt{a},\quad b,\quad c,\quad d,\quad e,\quad q^{-n} \\
		\multicolumn{2}{c}{
			\begin{matrix}
			   \sqrt{a},\,\, -\sqrt{a},\,\, aq/b,\,\, aq/c,\,\, aq/d,\,\, aq/e,\,\, aq^{n+1} 	
			\end{matrix}}
	\end{matrix} \hspace{8pt}
;q, q \right]\\ &\qquad\qquad = \dfrac{(aq,aq/bc,aq/bd,aq/cd;q)_n}{(aq/b,aq/c,aq/d,aq/bcd;q)_n},\numberthis\label{eq112}   
\end{align*}\\
where $a^2q = bcdeq^{-n}$.\\
\end{theorem}

\begin{theorem}(Andrews-Berkovich \cite[p. $535$, eq.$(3.1)$ and eq.$(3.2)$]{And-Ber02})\label{thm19}
For any non-negative integer $n$, we have\\
\begin{align*}
&{}_{10}\phi_9 \left[
	\setlength\arraycolsep{2pt}
	\begin{matrix}
		a, q\sqrt{a}, -q\sqrt{a}, a\sqrt{q/k}, -a\sqrt{q/k}, aq/\sqrt{k}, -aq/\sqrt{k}, k/aq, kq^n, q^{-n} \\
		\multicolumn{2}{c}{
			\begin{matrix}
			   \sqrt{a},\,\, -\sqrt{a},\,\, \sqrt{kq},\,\, -\sqrt{kq},\,\, \sqrt{k},\,\, -\sqrt{k},\,\, a^2q^2/k,\,\, aq^{-n+1}/k,\,\, aq^{n+1} 	
			\end{matrix}}
	\end{matrix} \hspace{2pt}
;q, q \right]\\ &\qquad\qquad = \dfrac{(aq,k^2/a^2q;q)_n}{(k,k/a;q)_n}\numberthis\label{eq113}   
\end{align*}\\
and\\
\begin{align*}
&{}_{10}\phi_9 \left[
	\setlength\arraycolsep{2pt}
	\begin{matrix}
		a, q\sqrt{a}, -q\sqrt{a}, a\sqrt{q/k}, -a\sqrt{q/k}, a/\sqrt{k}, -aq/\sqrt{k}, k/a, kq^n, q^{-n} \\
		\multicolumn{2}{c}{
			\begin{matrix}
			   \sqrt{a},\,\, -\sqrt{a},\,\, \sqrt{kq},\,\, -\sqrt{kq},\,\, q\sqrt{k},\,\, -\sqrt{k},\,\, a^2q/k,\,\, aq^{-n+1}/k,\,\, aq^{n+1} 	
			\end{matrix}}
	\end{matrix} \hspace{2pt}
;q, q \right]\\ &\qquad\qquad = \dfrac{(aq,\sqrt{k},k^2/a^2;q)_n}{(k,k/a,q\sqrt{k};q)_n}.\numberthis\label{eq114}
\end{align*}\\
\end{theorem}

\par In Section \ref{s2}, we prove Theorems \ref{thm11}, \ref{thm12}, \ref{thm13}, \ref{thm14}, \ref{thm16}, \ref{thm17}, \ref{thm18}, and \ref{thm19} using Pfaff's original method. We also show that upon applying a quadratic transformation, Theorem \ref{thm15} gives rise to Theorem \ref{thm13}.\\

\section{Proofs}\label{s2}
In this section, we present the proofs of Theorems \ref{thm11} - \ref{thm19}. From hereon in all the proofs, we want to remind the reader that while we will use $\infty$ as the upper limit of summation to save trouble in the re-indexing step, the sums are all finite.\\
\subsection{Proof of Theorem \ref{thm11}}\label{ss21}
Consider $S_n(z;q)$ to be the the left-hand side of \eqref{eq14}. Then, we have\\
\begin{align*}
&\quad S_n(z;q) - S_{n-1}(z;q)
\\
&= \sum\limits_{j = 0}^{\infty}\dfrac{(q^{-n};q)_j}{(q;q)_j}z^j - \sum\limits_{j = 0}^{\infty}\dfrac{(q^{-n+1};q)_j}{(q;q)_j}z^j
\\
&= \sum\limits_{j = 1}^{\infty}\dfrac{z^j}{(q;q)_j}\left((q^{-n};q)_j - (q^{-n+1};q)_j\right)
\\
&= -q^{-n}\sum\limits_{j = 1}^{\infty}\dfrac{(q^{-n+1};q)_{j-1}}{(q;q)_{j-1}}z^j
\\
&= -zq^{-n}\sum\limits_{j = 0}^{\infty}\dfrac{(q^{-n+1};q)_j}{(q;q)_j}z^j
\end{align*}
\newpage
\begin{align*}
\\
&= -zq^{-n}S_{n-1}(z;q).
\end{align*}\\
Thus, we have\\
\begin{align}\label{eq21}
S_n(z;q) = (1-zq^{-n})S_{n-1}(z;q).     
\end{align}\\
This already implies Theorem \ref{thm11} by iteration since $S_0(z;q) = 1$, but for the sake of completeness we continue with Pfaff’s method. Now, let $S_n^{\prime}(z;q)$ denote the right-hand side of \eqref{eq14}. Then, we have\\
\begin{align*}
&\quad S_n^{\prime}(z;q) - S_{n-1}^{\prime}(z;q)\\ &= (zq^{-n};q)_n - (zq^{-n+1};q)_{n-1}
\\
&= -zq^{-n}(zq^{-n+1};q)_{n-1}
\\
&= -zq^{-n}S_{n-1}^{\prime}(z;q).
\end{align*}\\
Thus, we have\\
\begin{align}\label{eq22}
S_n^{\prime}(z;q) = (1-zq^{-n})S_{n-1}^{\prime}(z;q).    
\end{align}\\
Now, observe that both $S_n(z;q)$ and $S_n^{\prime}(z;q)$ have the same initial conditions $S_0(z;q) = S_0^{\prime}(z;q) = 1$ and both obey the same recurrences \eqref{eq21} and \eqref{eq22} respectively. Hence,\\
\begin{align*}
S_n(z;q) = S_n^{\prime}(z;q)    
\end{align*}\\
and we have a proof of the $q$-binomial theorem (Theorem \ref{thm11}) using Pfaff's method.\qed\\

\subsection{Proof of Theorem \ref{thm12}}\label{ss22}
Consider $S_n(a,c;q)$ to be the the left-hand side of \eqref{eq15}. Then, we have\\
\begin{align*}
&\quad S_n(a,c;q) - S_{n-1}(a,c;q)
\\
&= \sum\limits_{j = 0}^{\infty}\dfrac{(a,q^{-n};q)_j}{(q,c;q)_j}q^j - \sum\limits_{j = 0}^{\infty}\dfrac{(a,q^{-n+1};q)_j}{(q,c;q)_j}q^j
\end{align*}
\newpage
\begin{align*}
\\
&= \sum\limits_{j = 1}^{\infty}\dfrac{(a;q)_jq^j}{(q,c;q)_j}\left((q^{-n};q)_j - (q^{-n+1};q)_j\right)
\\
&= -q^{-n}\sum\limits_{j = 1}^{\infty}\dfrac{(a;q)_j(q^{-n+1};q)_{j-1}}{(q;q)_{j-1}(c;q)_j}q^j
\\
&= -q^{-n+1}\sum\limits_{j = 0}^{\infty}\dfrac{(a;q)_{j+1}(q^{-n+1};q)_j}{(q;q)_j(c;q)_{j+1}}q^j
\\
&= -\dfrac{q^{-n+1}(1-a)}{(1-c)}\sum\limits_{j = 0}^{\infty}\dfrac{(aq,q^{-n+1};q)_j}{(q,cq;q)_j}q^j
\\
&= -\dfrac{q^{-n+1}(1-a)}{(1-c)}S_{n-1}(aq,cq;q).
\end{align*}\\
Thus, we have\\
\begin{align}\label{eq23}
S_n(a,c;q) - S_{n-1}(a,c;q) = -\dfrac{q^{-n+1}(1-a)}{(1-c)}S_{n-1}(aq,cq;q).     
\end{align}\\
Now, let $S_n^{\prime}(a,c;q)$ denote the right-hand side of \eqref{eq15}. Then, we have\\
\begin{align*}
&\quad S_n^{\prime}(a,c;q) - S_{n-1}^{\prime}(a,c;q)\\ &= \dfrac{(c/a;q)_n}{(c;q)_n}a^n - \dfrac{(c/a;q)_{n-1}}{(c;q)_{n-1}}a^{n-1}
\\
&= -\dfrac{(1-a)(c/a;q)_{n-1}}{(c;q)_n}a^{n-1}
\\
&= -\dfrac{q^{-n+1}(1-a)(cq/aq;q)_{n-1}}{(1-c)(cq;q)_{n-1}}(aq)^{n-1}
\\
&= -\dfrac{q^{-n+1}(1-a)}{(1-c)}S_{n-1}^{\prime}(aq,cq;q).
\end{align*}\\
Thus, we have\\
\begin{align}\label{eq24}
S_n^{\prime}(a,c;q) - S_{n-1}^{\prime}(a,c;q) = -\dfrac{q^{-n+1}(1-a)}{(1-c)}S_{n-1}^{\prime}(aq,cq;q).     
\end{align}\\
Now, observe that both $S_n(a,c;q)$ and $S_n^{\prime}(a,c;q)$ have the same initial conditions $S_0(a,c;q) = S_0^{\prime}(a,c;q) = 1$ and both obey the same recurrences \eqref{eq23} and \eqref{eq24} respectively. Hence,\\
\begin{align*}
S_n(a,c;q) = S_n^{\prime}(a,c;q)    
\end{align*}\\
and we have a proof of $q$-Chu-Vandermonde sum (Theorem \ref{thm12}) using Pfaff's method.\qed\\

\subsection{Proof of Theorem \ref{thm13}}\label{ss23}
Consider $S_n(a,b,c;q)$ to be the the left-hand side of \eqref{eq16}. Then, we have\\
\begin{align*}
&\quad S_n(a,b,c;q) - S_{n-1}(a,b,c;q)
\\ 
&= \sum\limits_{j = 0}^{\infty}\dfrac{(a,b,q^{-n};q)_j}{(q,c,abq^{-n+1}/c;q)_j}q^j - \dfrac{(a,b,q^{-n+1};q)j}{(q,c,abq^{-n+2}/c;q)_j}q^j
\\
&= \sum\limits_{j = 1}^{\infty}\dfrac{(a,b;q)_jq^j}{(q,c;q)_j}\left(\dfrac{(q^{-n};q)_j}{(abq^{-n+1}/c;q)_j} - \dfrac{(q^{-n+1};q)_j}{(abq^{-n+2}/c;q)_j}\right)
\\
&= -\dfrac{q^{-n}(1-abq/c)}{(1-abq^{-n+1}/c)}\sum\limits_{j = 1}^{\infty}\dfrac{(a,b;q)_j(q^{-n+1};q)_{j-1}}{(q;q)_{j-1}(c,abq^{-n+2}/c;q)_j}q^j
\\
&= -\dfrac{q^{-n}(1-abq/c)}{(1-abq^{-n+1}/c)}\sum\limits_{j = 0}^{\infty}\dfrac{(a,b;q)_{j+1}(q^{-n+1};q)_j}{(q;q)_j(c,abq^{-n+2}/c;q)_{j+1}}q^{j+1}
\\
&= -\dfrac{q^{-n+1}(1-a)(1-b)(1-abq/c)}{(1-c)(1-abq^{-n+1}/c)(1-abq^{-n+2}/c)}\\ &\qquad\qquad\times\sum\limits_{j = 0}^{\infty}\dfrac{(aq,bq,q^{-n+1};q)_j}{(q,cq,abq^{-n+3}/c;q)_j}q^j
\\
&= -\dfrac{q^{-n+1}(1-a)(1-b)(1-abq/c)}{(1-c)(1-abq^{-n+1}/c)(1-abq^{-n+2}/c)}S_{n-1}(aq,bq,cq;q).
\end{align*}\\
Thus, we have\\
\begin{align*}
&\quad S_n(a,b,c;q) - S_{n-1}(a,b,c;q)\\ &= -\dfrac{q^{-n+1}(1-a)(1-b)(1-abq/c)}{(1-c)(1-abq^{-n+1}/c)(1-abq^{-n+2}/c)}S_{n-1}(aq,bq,cq;q).\numberthis\label{eq25}    
\end{align*}\\
\begin{remark}
\eqref{eq25} was also proved by Andrews \cite[p. $22$, eq.$(7.3)$]{And96}.\\
\end{remark}
Now, let $S_n^{\prime}(a,b,c;q)$ denote the right-hand side of \eqref{eq16}. Then, we have\\
\begin{align*}
&\quad S_n^{\prime}(a,b,c;q) - S_{n-1}^{\prime}(a,b,c;q)\\ &= \dfrac{(c/a,c/b;q)_n}{(c,c/ab;q)_n} - \dfrac{(c/a,c/b;q)_{n-1}}{(c,c/ab;q)_{n-1}}
\\
&= \dfrac{cq^{n-1}(1-1/a)(1-1/b)(c/a,c/b;q)_{n-1}}{(c,c/ab;q)_n}
\\
&= \dfrac{cq^{n-1}(1-a)(1-b)(c/abq;q)_{n-1}(cq/aq,cq/bq;q)_{n-1}}{ab(1-c)(1-c/ab)(cq/ab;q)_{n-1}(cq,c/abq;q)_{n-1}}
\\
&= \dfrac{cq^{n-1}(1-a)(1-b)(1-c/abq)}{ab(1-c)(1-cq^{n-2}/ab)(1-cq^{n-1}/ab)}S_{n-1}^{\prime}(aq,bq,cq;q)
\\
&= -\dfrac{q^{-n+1}(1-a)(1-b)(1-abq/c)}{(1-c)(1-abq^{-n+1}/c)(1-abq^{-n+2}/c)}S_{n-1}^{\prime}(aq,bq,cq;q).
\end{align*}\\
Thus, we have\\
\begin{align*}
&\quad S_n^{\prime}(a,b,c;q) - S_{n-1}^{\prime}(a,b,c;q)\\ &= -\dfrac{q^{-n+1}(1-a)(1-b)(1-abq/c)}{(1-c)(1-abq^{-n+1}/c)(1-abq^{-n+2}/c)}S_{n-1}^{\prime}(aq,bq,cq;q).\numberthis\label{eq26}    
\end{align*}\\
Now, observe that both $S_n(a,b,c;q)$ and $S_n^{\prime}(a,b,c;q)$ have the same initial conditions $S_0(a,b,c;q)$ $=$ $S_0^{\prime}(a,b,c;q) = 1$ and both obey the same recurrences \eqref{eq25} and \eqref{eq26} respectively. Hence,\\
\begin{align*}
S_n(a,b,c;q) = S_n^{\prime}(a,b,c;q)    
\end{align*}\\
and we have a proof of $q$-Pfaff-Saalsch\"utz sum (Theorem \ref{thm13}) using Pfaff's method.\qed\\

\subsection{Proof of Theorem \ref{thm14}}\label{ss24}
Consider $S_n(a,b;q)$ to be the the left-hand side of \eqref{eq17}. Then, we have\\
\begin{align*}
&\quad S_n(a,b;q) - S_{n-1}(a,b;q)
\\
&= \sum\limits_{j = 0}^{\infty}\dfrac{(a,-q\sqrt{a},b,q^{-n};q)_j}{(q,-\sqrt{a},aq/b,aq^{n+1};q)_j}\left(\dfrac{q^{n+1}\sqrt{a}}{b}\right)^j\\ &\qquad\qquad - \sum\limits_{j = 0}^{\infty}\dfrac{(a,-q\sqrt{a},b,q^{-n+1};q)_j}{(q,-\sqrt{a},aq/b,aq^n;q)_j}\left(\dfrac{q^n\sqrt{a}}{b}\right)^j
\end{align*}
\newpage
\begin{align*}
\\
&= \sum\limits_{j = 1}^{\infty}\dfrac{(a,-q\sqrt{a},b;q)_j(q^n\sqrt{a}/b)^j}{(q,-\sqrt{a},aq/b;q)_j}\left(\dfrac{(q^{-n};q)_jq^j}{(aq^{n+1};q)_j}-\dfrac{(q^{-n+1};q)_j}{(aq^n;q)_j}\right)
\\
&= -\dfrac{1}{(1-aq^n)}\sum\limits_{j = 1}^{\infty}\dfrac{(a,-q\sqrt{a},b;q)_j(q^{-n+1};q)_{j-1}(1-aq^j)}{(-\sqrt{a},aq/b;q)_j(q,aq^{n+1};q)_{j-1}(1-aq^{n+j})}\\ &\qquad\qquad\times\left(\dfrac{q^n\sqrt{a}}{b}\right)^j
\\
&= -\dfrac{(1-a)}{(1-aq^n)}\sum\limits_{j = 1}^{\infty}\dfrac{(aq,-q\sqrt{a},b;q)_j(q^{-n+1};q)_{j-1}}{(q;q)_{j-1}(-\sqrt{a},aq/b,aq^{n+1};q)_j}\left(\dfrac{q^n\sqrt{a}}{b}\right)^j
\\
&= -\dfrac{(1-a)}{(1-aq^n)}\sum\limits_{j = 0}^{\infty}\dfrac{(aq,-q\sqrt{a},b;q)_{j+1}(q^{-n+1};q)_j}{(q;q)_j(-\sqrt{a},aq/b,aq^{n+1};q)_{j+1}}\left(\dfrac{q^n\sqrt{a}}{b}\right)^{j+1}
\\
&= -\dfrac{(1-a)(1-aq)(1+q\sqrt{a})(1-b)(q^n\sqrt{a}/b)}{(1-aq^n)(1+\sqrt{a})(1-aq/b)(1-aq^{n+1})}\\ &\qquad\qquad\times\sum\limits_{j = 0}^{\infty}\dfrac{(aq^2,-q^2\sqrt{a},bq,q^{-n+1};q)_j}{(q,-q\sqrt{a},aq^2/b,aq^{n+2};q)_j}\left(\dfrac{q^n\sqrt{a}}{b}\right)^j
\\
&= -\dfrac{(1-\sqrt{a})(1+q\sqrt{a})(1-aq)(1-b)(q^n\sqrt{a}/b)}{(1-aq/b)(1-aq^n)(1-aq^{n+1})}S_{n-1}(aq^2,bq;q).
\end{align*}\\
Thus, we have\\
\begin{align*}
&\quad S_n(a,b;q) - S_{n-1}(a,b;q)\\ &= -\dfrac{(1-\sqrt{a})(1+q\sqrt{a})(1-aq)(1-b)(q^n\sqrt{a}/b)}{(1-aq/b)(1-aq^n)(1-aq^{n+1})}S_{n-1}(aq^2,bq;q).\numberthis\label{eq27}    
\end{align*}\\
Now, let $S_n^{\prime}(a,b;q)$ denote the right-hand side of \eqref{eq17}. Then, we have\\
\begin{align*}
&\quad S_n^{\prime}(a,b;q) - S_{n-1}^{\prime}(a,b;q)
\\
&= \dfrac{(aq,q\sqrt{a}/b;q)_n}{(q\sqrt{a},aq/b;q)_n} - \dfrac{(aq,q\sqrt{a}/b;q)_{n-1}}{(q\sqrt{a},aq/b;q)_{n-1}}
\\
&= \dfrac{q^n\sqrt{a}(1-\sqrt{a})(1-1/b)(aq,q\sqrt{a}/b;q)_{n-1}}{(q\sqrt{a},aq/b;q)_n}
\\
&= -\dfrac{(1-\sqrt{a})(1-b)(q^n\sqrt{a}/b)(aq,q\sqrt{a}/b,q^2\sqrt{a},aq^2/b;q)_{n-1}}{(q\sqrt{a},aq/b;q)_{n}(aq^3,q\sqrt{a}/b;q)_{n-1}}\\&\qquad\qquad\times\dfrac{(aq^3,q\sqrt{a}/b;q)_{n-1}}{(q^2\sqrt{a},aq^2/b;q)_{n-1}}
\end{align*}
\newpage
\begin{align*}
\\
&= -\dfrac{(1-\sqrt{a})(1-b)(q^n\sqrt{a}/b)(aq;q)_{n-1}}{(1-q\sqrt{a})(1-aq/b)(aq^3;q)_{n-1}}S_{n-1}^{\prime}(aq^2,bq;q)
\\
&= -\dfrac{(1-\sqrt{a})(1+q\sqrt{a})(1-aq)(1-b)(q^n\sqrt{a}/b)}{(1-aq/b)(1-aq^n)(1-aq^{n+1})}S_{n-1}^{\prime}(aq^2,bq;q).
\end{align*}\\
Thus, we have\\
\begin{align*}
&\quad S_n^{\prime}(a,b;q) - S_{n-1}^{\prime}(a,b;q)\\ &= -\dfrac{(1-\sqrt{a})(1+q\sqrt{a})(1-aq)(1-b)(q^n\sqrt{a}/b)}{(1-aq/b)(1-aq^n)(1-aq^{n+1})}S_{n-1}^{\prime}(aq^2,bq;q).\numberthis\label{eq28}
\end{align*}\\
Now, observe that both $S_n(a,b;q)$ and $S_n^{\prime}(a,b;q)$ have the same initial conditions $S_0(a,b;q)$ $=$ $S_0^{\prime}(a,b;q) = 1$ and both obey the same recurrences \eqref{eq27} and \eqref{eq28} respectively. Hence,\\
\begin{align*}
S_n(a,b;q) = S_n^{\prime}(a,b;q)    
\end{align*}\\
and we have a proof of $q$-Dixon sum (Theorem \ref{thm14}) using Pfaff's method.\qed\\

\subsection{Proof of Theorem \ref{thm15}}\label{ss25}
Consider $S_n(a,c;q)$ to be the the left-hand side of \eqref{eq18}. Then, using the substitution $(a,b,c,d)\longmapsto (\sqrt{a},\sqrt{aq},-c,q^{-n})$ in Singh's quadratic transformation \cite[p. 361, III.$21$]{Gas-Rah04} (change of base from $q^2\rightarrow q$) on $S_n(a,c;q)$, we get\\
\begin{align*}
&\quad S_n(a,c;q)\\ &= {}_4\phi_3 \left[
	\setlength\arraycolsep{2pt}
	\begin{matrix}
		a,\quad aq,\quad c^2,\quad q^{-2n} \\
		\multicolumn{2}{c}{
			\begin{matrix}
			  \,\,\,\,\, a^2q^2,\,\,\, cq^{-n},\,\,\, cq^{-n+1} 	
			\end{matrix}}
	\end{matrix} \hspace{8pt}
;q^2, q^2 \right]
\\
&= {}_4\phi_3 \left[
	\setlength\arraycolsep{2pt}
	\begin{matrix}
		a,\quad aq,\quad -c,\quad q^{-n} \\
		\multicolumn{2}{c}{
			\begin{matrix}
			  \,\,\,\,\, aq,\,\,\, -aq,\,\,\, cq^{-n} 	
			\end{matrix}}
	\end{matrix} \hspace{12pt}
;q, q \right]
\\
&= {}_3\phi_2 \left[
	\setlength\arraycolsep{2pt}
	\begin{matrix}
		a,\,\,\,\,\, -c,\,\,\,\,\, q^{-n} \\
		\multicolumn{2}{c}{
			\begin{matrix}
			  \,\,\,\,\, -aq,\,\,\, cq^{-n} 	
			\end{matrix}}
	\end{matrix} \hspace{8pt}
;q, q \right]
\\
&= \dfrac{(-q,aq/c;q)_n}{(-aq,q/c;q)_n}
\end{align*}\\
where the last line follows by substituting $(a,b,c)\longmapsto (a,-c,-aq)$ in \eqref{eq16}. Thus, we obtain a direct proof of Theorem \ref{thm15} without resorting to Pfaff's method.\qed\\

\subsection{Proof of Theorem \ref{thm16}}\label{ss26}
Consider $A_n(b,d;q)$ and $B_n(b,d;q)$ to be the left-hand sides of \eqref{eq19} and \eqref{eq110} respectively. Then, we have\\
\begin{align*}
&\quad B_n(b,d;q) - A_n(b,d;q)\\ &= \sum\limits_{j = 0}^{\infty}\dfrac{(q^{-n},b,b\sqrt{q},d^2q^{n+1};q)_j}{(q,dq,d\sqrt{q},b^2q;q)_j}q^j - \sum\limits_{j = 0}^{\infty}\dfrac{(q^{-n},b,b\sqrt{q},d^2q^n;q)_j}{(q,dq,d\sqrt{q},b^2;q)_j}q^j
\\
&= \sum\limits_{j = 1}^{\infty}\dfrac{(q^{-n},b,b\sqrt{q};q)_jq^j}{(q,dq,d\sqrt{q};q)_j}\left(\dfrac{(d^2q^{n+1};q)_j}{(b^2q;q)_j}-\dfrac{(d^2q^n;q)_j}{(b^2;q)_j}\right)
\\
&= -\dfrac{(b^2-d^2q^n)}{(1-b^2)}\sum\limits_{j = 1}^{\infty}\dfrac{(q^{-n},b,b\sqrt{q};q)_j(d^2q^{n+1};q)_{j-1}}{(q;q)_{j-1}(dq,d\sqrt{q},b^2q;q)_j}q^j
\\
&= -\dfrac{(b^2-d^2q^n)}{(1-b^2)}\sum\limits_{j = 0}^{\infty}\dfrac{(q^{-n},b,b\sqrt{q};q)_{j+1}(d^2q^{n+1};q)_j}{(q;q)_j(dq,d\sqrt{q},b^2q;q)_{j+1}}q^{j+1}
\\
&= -\dfrac{q(b^2-d^2q^n)(1-q^{-n})(1-b)(1-b\sqrt{q})}{(1-b^2)(1-dq)(1-d\sqrt{q})(1-b^2q)}\\&\qquad\qquad\times\sum\limits_{j = 0}^{\infty}\dfrac{(q^{-n+1},bq,bq\sqrt{q},d^2q^{n+1};q)_j}{(q,dq^2,dq\sqrt{q},b^2q^2;q)_j}q^j
\\
&= -\dfrac{q(b^2-d^2q^n)(1-q^{-n})}{(1+b)(1+b\sqrt{q})(1-dq)(1-d\sqrt{q})}A_{n-1}(bq,dq;q).
\end{align*}\\
Thus, we have\\
\begin{align*}
&\quad B_n(b,d;q) - A_n(b,d;q)\\ &= -\dfrac{q(b^2-d^2q^n)(1-q^{-n})}{(1+b)(1+b\sqrt{q})(1-dq)(1-d\sqrt{q})}A_{n-1}(bq,dq;q),\numberthis\label{eq29}   
\end{align*}\\
while it can similarly be shown that\\
\begin{align*}
&\quad A_n(b,d;q) - B_{n-1}(b,d;q)\\ &= \dfrac{q(1-d^2q^n)(b^2-q^{-n})}{(1+b)(1+b\sqrt{q})(1-dq)(1-d\sqrt{q})}A_{n-1}(bq,dq;q).\numberthis\label{eq210}
\end{align*}\\
Now, let $A_n^{\prime}(b,d;q)$ and $B_n^{\prime}(b,d;q)$ denote the right-hand sides of \eqref{eq19} and \eqref{eq110} respectively. Then, we have\\
\begin{align*}
&\quad B_n^{\prime}(b,d;q) - A_n^{\prime}(b,d;q)\\ &= \dfrac{b^n(-\sqrt{q},d\sqrt{q}/b;\sqrt{q})_n}{(-b\sqrt{q},d\sqrt{q};\sqrt{q})_n} - \dfrac{b^n(1-d)(-\sqrt{q},d\sqrt{q}/b;\sqrt{q})_n}{(1-q^nd)(-b,d;\sqrt{q})_n}
\\
&= \dfrac{b^n(-\sqrt{q},d\sqrt{q}/b;\sqrt{q})_n}{(-b\sqrt{q},d\sqrt{q};\sqrt{q})_{n-1}}\left(\dfrac{1}{(1+bq^{n/2})(1-dq^{n/2})}-\dfrac{1}{(1-q^nd)(1+b)}\right)
\\
&= -\dfrac{q(b^2-d^2q^n)(1-q^{-n})}{(1+b)(1+b\sqrt{q})(1-dq)(1-d\sqrt{q})}\\&\qquad\qquad\times\dfrac{(bq)^{n-1}(1-dq)(-\sqrt{q},d\sqrt{q}/b;\sqrt{q})_{n-1}}{(1-dq^n)(-bq,dq;\sqrt{q})_{n-1}}
\\
&= -\dfrac{q(b^2-d^2q^n)(1-q^{-n})}{(1+b)(1+b\sqrt{q})(1-dq)(1-d\sqrt{q})}A_{n-1}^{\prime}(bq,dq;q).
\end{align*}\\
Thus, we have\\
\begin{align*}
&\quad B_n^{\prime}(b,d;q) - A_n^{\prime}(b,d;q)\\ &= -\dfrac{q(b^2-d^2q^n)(1-q^{-n})}{(1+b)(1+b\sqrt{q})(1-dq)(1-d\sqrt{q})}A_{n-1}^{\prime}(bq,dq;q),\numberthis\label{eq211}
\end{align*}\\
while it can similarly be shown that\\
\begin{align*}
&\quad A_n^{\prime}(b,d;q) - B_{n-1}^{\prime}(b,d;q)\\ &= \dfrac{q(1-d^2q^n)(b^2-q^{-n})}{(1+b)(1+b\sqrt{q})(1-dq)(1-d\sqrt{q})}A_{n-1}^{\prime}(bq,dq;q).\numberthis\label{eq212}   
\end{align*}\\
Now, observe that \eqref{eq29} and \eqref{eq210} together with the initial conditions $A_0(b,d;q)$ $=$ $B_0(b,d;q) = 1$ uniquely define both $A_n(b,d;q)$ and $B_n(b,d;q)$. Similarly, \eqref{eq211} and \eqref{eq212} together with the initial conditions $A_0^{\prime}(b,d;q)$ $=$ $B_0^{\prime}(b,d;q) = 1$ uniquely define both $A_n^{\prime}(b,d;q)$ and $B_n^{\prime}(b,d;q)$. Also, note that the recurrences \eqref{eq29} and \eqref{eq211} are the same and the recurrences \eqref{eq210} and \eqref{eq212} are the same. Hence,\\
\begin{align*}
A_n(b,d;q) &= A_n^{\prime}(b,d;q)
\end{align*}
and
\begin{align*}
B_n(b,d;q) &= B_n^{\prime}(b,d;q).
\end{align*}\\
Thus, we have a proof of the two ${}_4\phi_3$ identities due to Andrews and Li-Chu (Theorem \ref{thm16}) using Pfaff's method.\qed\\

\subsection{Proof of Theorem \ref{thm17}}\label{ss27}
Consider $S_n(a,b,c;q)$ to be the left-hand side of \eqref{eq111}. Then, we have\\
\begin{align*}
&\quad S_n(a,b,c;q) - S_{n-1}(a,b,c;q)\\ &= \sum\limits_{j = 0}^{\infty}\dfrac{(a,q\sqrt{a},-q\sqrt{a},b,c,q^{-n};q)_j}{(q,\sqrt{a},-\sqrt{a},aq/b,aq/c,aq^{n+1};q)_j}\left(\dfrac{aq^{n+1}}{bc}\right)^j\\&\qquad\qquad - \sum\limits_{j = 0}^{\infty}\dfrac{(a,q\sqrt{a},-q\sqrt{a},b,c,q^{-n+1};q)_j}{(q,\sqrt{a},-\sqrt{a},aq/b,aq/c,aq^n;q)_j}\left(\dfrac{aq^n}{bc}\right)^j
\\
&= \sum\limits_{j = 1}^{\infty}\dfrac{(a,q\sqrt{a},-q\sqrt{a},b,c;q)_j(aq^n/bc)^j}{(q,\sqrt{a},-\sqrt{a},aq/b,aq/c;q)_j}\\ &\qquad\qquad\times\left(\dfrac{(q^{-n};q)_jq^j}{(aq^{n+1};q)_j} - \dfrac{(q^{-n+1};q)_j}{(aq^n;q)_j}\right)
\\
&= -\dfrac{(1-a)}{(1-aq^n)}\sum\limits_{j = 1}^{\infty}\dfrac{(aq,q\sqrt{a},-q\sqrt{a},b,c;q)_j(q^{-n+1};q)_{j-1}}{(q;q)_{j-1}(\sqrt{a},-\sqrt{a},aq/b,aq/c,aq^{n+1};q)_j}\\ &\qquad\qquad\times\left(\dfrac{aq^n}{bc}\right)^j
\\
&= -\dfrac{(1-a)}{(1-aq^n)}\sum\limits_{j = 0}^{\infty}\dfrac{(aq,q\sqrt{a},-q\sqrt{a},b,c;q)_{j+1}(q^{-n+1};q)_j}{(q;q)_j(\sqrt{a},-\sqrt{a},aq/b,aq/c,aq^{n+1};q)_{j+1}}\\ &\qquad\qquad\times\left(\dfrac{aq^n}{bc}\right)^{j+1}
\\
&= -\dfrac{(1-a)(1-aq)(1-q\sqrt{a})(1+q\sqrt{a})(1-b)(1-c)(aq^n/bc)}{(1-\sqrt{a})(1+\sqrt{a})(1-aq/b)(1-aq/c)(1-aq^n)(1-aq^{n+1})}\\&\qquad\qquad\times\sum\limits_{j = 0}^{\infty}\dfrac{(aq^2,q^2\sqrt{a},-q^2\sqrt{a},bq,cq,q^{-n+1};q)_j}{(q,q\sqrt{a},-q\sqrt{a},aq^2/b,aq^2/c,aq^{n+2};q)_j}\left(\dfrac{aq^n}{bc}\right)^j
\\
&= -\dfrac{(1-aq)(1-aq^2)(1-b)(1-c)(aq^n/bc)}{(1-aq^n)(1-aq^{n+1})(1-aq/b)(1-aq/c)}S_{n-1}(aq^2,bq,cq;q).
\end{align*}\\
Thus, we have\\
\begin{align*}
&\quad S_n(a,b,c;q) - S_{n-1}(a,b,c;q)\\ &= -\dfrac{(1-aq)(1-aq^2)(1-b)(1-c)(aq^n/bc)}{(1-aq^n)(1-aq^{n+1})(1-aq/b)(1-aq/c)}S_{n-1}(aq^2,bq,cq;q).\numberthis\label{eq213}
\end{align*}\\
Now, let $S_n^{\prime}(a,b,c;q)$ denote the right-hand side of \eqref{eq111}. Then, we have\\
\begin{align*}
&\quad S_{n}^{\prime}(a,b,c;q) - S_{n-1}^{\prime}(a,b,c;q)\\ &= \dfrac{(aq,aq/bc;q)_n}{(aq/b,aq/c;q)_n} - \dfrac{(aq,aq/bc;q)_{n-1}}{(aq/b,aq/c;q)_{n-1}}
\\
&= -\dfrac{aq^n(1-1/b)(1-1/c)(aq,aq/bc;q)_{n-1}}{(aq/b,aq/c;q)_n}
\\
&= -\dfrac{(1-b)(1-c)(aq^n/bc)(aq,aq/bc,aq^2/b,aq^2/c;q)_{n-1}}{(aq/b,aq/c;q)_n(aq^3,aq/bc;q)_{n-1}}\\&\qquad\qquad\times\dfrac{(aq^3,aq/bc;q)_{n-1}}{(aq^2/b,aq^2/c;q)_{n-1}}
\\
&= -\dfrac{(1-b)(1-c)(aq^n/bc)(aq,aq/bc,aq^2/b,aq^2/c;q)_{n-1}}{(aq/b,aq/c;q)_n(aq^3,aq/bc;q)_{n-1}}\\&\qquad\qquad\times S_{n-1}^{\prime}(aq^2,bq,cq;q)
\\
&= -\dfrac{(1-aq)(1-aq^2)(1-b)(1-c)(aq^n/bc)}{(1-aq^n)(1-aq^{n+1})(1-aq/b)(1-aq/c)}S_{n-1}^{\prime}(aq^2,bq,cq;q).
\end{align*}\\
Thus, we have\\
\begin{align*}
&\quad S_n^{\prime}(a,b,c;q) - S_{n-1}^{\prime}(a,b,c;q)\\ &= -\dfrac{(1-aq)(1-aq^2)(1-b)(1-c)(aq^n/bc)}{(1-aq^n)(1-aq^{n+1})(1-aq/b)(1-aq/c)}S_{n-1}^{\prime}(aq^2,bq,cq;q).\numberthis\label{eq214}
\end{align*}\\
Now, observe that both $S_n(a,b,c;q)$ and $S_n^{\prime}(a,b,c;q)$ have the same initial conditions $S_0(a,b,c;q) = S_0^{\prime}(a,b,c;q) = 1$ and both obey the same recurrences \eqref{eq213} and \eqref{eq214} respectively. Hence,\\
\begin{align*}
S_n(a,b,c;q) = S_n^{\prime}(a,b,c;q)    
\end{align*}\\
and we have a proof of Rogers' ${}_6\phi_5$ sum (Theorem \ref{thm17}) using Pfaff's method.\qed\\

\subsection{Proof of Theorem \ref{thm18}}\label{ss28}
Consider $S_n(a,b,c,d;q)$ to be the left-hand side of \eqref{eq112}. Then, we have\\
\begin{align*}
&\quad S_n(a,b,c,d;q) - S_{n-1}(a,b,c,d;q)\\ &= \sum\limits_{j = 0}^{\infty}\dfrac{(a,q\sqrt{a},-q\sqrt{a},b,c,d,a^2q^{n+1}/bcd,q^{-n};q)_j}{(q,\sqrt{a},-\sqrt{a},aq/b,aq/c,aq/d,bcdq^{-n}/a,aq^{n+1};q)_j}q^j\\&\qquad - \sum\limits_{j = 0}^{\infty}\dfrac{(a,q\sqrt{a},-q\sqrt{a},b,c,d,a^2q^n/bcd,q^{-n+1};q)_j}{(q,\sqrt{a},-\sqrt{a},aq/b,aq/c,aq/d,bcdq^{-n+1}/a,aq^n;q)_j}q^j
\\
&= \sum\limits_{j = 1}^{\infty}\dfrac{(a,q\sqrt{a},-q\sqrt{a},b,c,d;q)_jq^j}{(q,\sqrt{a},-\sqrt{a},aq/b,aq/c,aq/d;q)_j}\\&\qquad\times\left(\dfrac{(a^2q^{n+1}/bcd,q^{-n};q)_j}{(bcdq^{-n}/a,aq^{n+1};q)_j} - \dfrac{(a^2q^n/bcd,q^{-n+1};q)_j}{(bcdq^{-n+1}/a,aq^n;q)_j}\right)
\\
&= -\dfrac{(q^{-n}-a^2q^n/bcd+aq^n-bcdq^{-n}/a)}{(1-bcdq^{-n}/a)(1-aq^n)}\\&\qquad\times\sum\limits_{j = 1}^{\infty}\dfrac{(a,q\sqrt{a},-q\sqrt{a},b,c,d;q)_j(a^2q^{n+1}/bcd,q^{-n+1};q)_{j-1}(1-aq^j)}{(q;q)_{j-1}(\sqrt{a},-\sqrt{a},aq/b,aq/c,aq/d,bcdq^{-n+1}/a,aq^{n+1};q)_j}q^j
\\
&= -\dfrac{(1-a)(q^{-n}-a^2q^n/bcd+aq^n-bcdq^{-n}/a)}{(1-bcdq^{-n}/a)(1-aq^n)}\\&\qquad\times\sum\limits_{j = 1}^{\infty}\dfrac{(aq,q\sqrt{a},-q\sqrt{a},b,c,d;q)_j(a^2q^{n+1}/bcd,q^{-n+1};q)_{j-1}}{(q;q)_{j-1}(\sqrt{a},-\sqrt{a},aq/b,aq/c,aq/d,bcdq^{-n+1}/a,aq^{n+1};q)_j}q^j
\\
&= -\dfrac{(1-a)(q^{-n}-a^2q^n/bcd+aq^n-bcdq^{-n}/a)}{(1-bcdq^{-n}/a)(1-aq^n)}\\&\qquad\times\sum\limits_{j = 0}^{\infty}\dfrac{(aq,q\sqrt{a},-q\sqrt{a},b,c,d;q)_{j+1}(a^2q^{n+1}/bcd,q^{-n+1};q)_j}{(q;q)_j(\sqrt{a},-\sqrt{a},aq/b,aq/c,aq/d,bcdq^{-n+1}/a,aq^{n+1};q)_{j+1}}q^{j+1}
\\
&= -\dfrac{q(q^{-n}-a^2q^n/bcd+aq^n-bcdq^{-n}/a)}{(1-bcdq^{-n}/a)(1-aq^n)}\\&\qquad\times\dfrac{(1-aq)(1-aq^2)(1-b)(1-c)(1-d)}{(1-aq/b)(1-aq/c)(1-aq/d)(1-bcdq^{-n+1}/a)(1-aq^{n+1})}\\&\qquad\times\sum\limits_{j = 0}^{\infty}\dfrac{(aq^2,q^2\sqrt{a},-q^2\sqrt{a},bq,cq,dq,a^2q^{n+1}/bcd,q^{-n+1};q)_j}{(q,q\sqrt{a},-q\sqrt{a},aq^2/b,aq^2/c,aq^2/d,bcdq^{-n+2}/a,aq^{n+2};q)_j}q^j
\\
&= -\dfrac{q(q^{-n}-a^2q^n/bcd+aq^n-bcdq^{-n}/a)}{(1-bcdq^{-n}/a)(1-aq^n)}\\&\qquad\times\dfrac{(1-aq)(1-aq^2)(1-b)(1-c)(1-d)}{(1-aq/b)(1-aq/c)(1-aq/d)(1-bcdq^{-n+1}/a)(1-aq^{n+1})}\\&\qquad\times S_{n-1}(aq^2,bq,cq,dq;q).
\end{align*}\\
Thus, we have\\
\begin{align}\label{eq215}
S_n(a,b,c,d;q) - S_{n-1}(a,b,c,d;q) = \phi_n(a,b,c,d;q)S_{n-1}(aq^2,bq,cq,dq;q),
\end{align}\\
where\\
\begin{align*}
&\quad \phi_n(a,b,c,d;q)\\ &= -\dfrac{q(1-aq)(1-aq^2)(1-b)(1-c)(1-d)}{(1-aq^n)(1-aq^{n+1})(1-aq/b)(1-aq/c)(1-aq/d)}\\ &\qquad\qquad\times \dfrac{(q^{-n}-a^2q^n/bcd+aq^n-bcdq^{-n}/a)}{(1-bcdq^{-n}/a)(1-bcdq^{-n+1}/a)}.    
\end{align*}\\
Now, let $S_n^{\prime}(a,b,c,d;q)$ denote the right-hand side of \eqref{eq112}. Then, we have\\
\begin{align*}
&\quad S_n^{\prime}(a,b,c,d;q) - S_{n-1}^{\prime}(a,b,c,d;q)\\ &= \dfrac{(aq,aq/bc,aq/bd,aq/cd;q)_n}{(aq/b,aq/c,aq/d,aq/bcd;q)_n} - \dfrac{(aq,aq/bc,aq/bd,aq/cd;q)_{n-1}}{(aq/b,aq/c,aq/d,aq/bcd;q)_{n-1}}
\\
&= -aq^n(1-a^2q^{2n}/bcd)(1-1/b)(1-1/c)(1-1/d)\\&\qquad\times\dfrac{(aq,aq/bc,aq/bd,aq/cd;q)_{n-1}}{(aq/b,aq/c,aq/d,aq/bcd;q)_n}
\\
&= -aq^n(1-a^2q^{2n}/bcd)(1-1/b)(1-1/c)(1-1/d)\\&\qquad\times\dfrac{(aq^2/b,aq^2/c,aq^2/d,a/bcd,aq,aq/bc,aq/bd,aq/cd;q)_{n-1}}{(aq^3,aq/bc,aq/bd,aq/cd;q)_{n-1}(aq/b,aq/c,aq/d,aq/bcd;q)_n}\\&\qquad\times\dfrac{(aq^3,aq/bc,aq/bd,aq/cd;q)_{n-1}}{(aq^2/b,aq^2/c,aq^2/d,a/bcd;q)_{n-1}}
\\
&= -\dfrac{aq^n(1-a^2q^{2n}/bcd)(1-1/b)(1-1/c)(1-1/d)(aq,a/bcd;q)_{n-1}}{(1-aq/b)(1-aq/c)(1-aq/d)(1-aq/bcd)(aq^3,aq^2/bcd;q)_{n-1}}\\&\qquad\times S_{n-1}^{\prime}(aq^2,bq,cq,dq;q)
\\
&= -\dfrac{q(q^{-n}-a^2q^n/bcd+aq^n-bcdq^{-n}/a)}{(1-bcdq^{-n}/a)(1-aq^n)}\\&\qquad\times\dfrac{(1-aq)(1-aq^2)(1-b)(1-c)(1-d)}{(1-aq/b)(1-aq/c)(1-aq/d)(1-bcdq^{-n+1}/a)(1-aq^{n+1})}\\&\qquad\times S_{n-1}^{\prime}(aq^2,bq,cq,dq;q).
\end{align*}\\
Thus, we have\\
\begin{align}\label{eq216}
S_n^{\prime}(a,b,c,d;q) - S_{n-1}^{\prime}(a,b,c,d;q) = \phi_n^{\prime}(a,b,c,d;q)S_{n-1}^{\prime}(aq^2,bq,cq,dq;q),
\end{align}\\
where\\
\begin{align*}
&\quad \phi_n^{\prime}(a,b,c,d;q)\\ &= -\dfrac{q(1-aq)(1-aq^2)(1-b)(1-c)(1-d)}{(1-aq^n)(1-aq^{n+1})(1-aq/b)(1-aq/c)(1-aq/d)}\\ &\qquad\qquad\times\dfrac{(q^{-n}-a^2q^n/bcd+aq^n-bcdq^{-n}/a)}{(1-bcdq^{-n}/a)(1-bcdq^{-n+1}/a)}.    
\end{align*}\\
Now, observe that both $S_n(a,b,c,d;q)$ and $S_n^{\prime}(a,b,c,d;q)$ have the same initial conditions $S_0(a,b,c,d;q)$ $=$\\ $S_0^{\prime}(a,b,c,d;q) = 1$ and both obey the same recurrences \eqref{eq215} and \eqref{eq216} respectively since $\phi_n(a,b,c,d;q) = \phi_n^{\prime}(a,b,c,d;q)$. Hence,\\
\begin{align*}
S_n(a,b,c,d;q) = S_n^{\prime}(a,b,c,d;q)    
\end{align*}\\
and we have a proof of Jackson's ${}_8\phi_7$ sum (Theorem \ref{thm18}) using Pfaff's method.\qed\\

\subsection{Proof of Theorem \ref{thm19}}\label{ss29}
Consider $A_n(a,k;q)$ and $B_n(a,k;q)$ to be the left-hand sides of \eqref{eq113} and \eqref{eq114} respectively. Then, we have\\
\begin{align*}
&\quad A_n(a,k;q) - B_n(a,k;q)
\\
&= \sum\limits_{j = 0}^{\infty}\dfrac{(a,q\sqrt{a},-q\sqrt{a},a\sqrt{q/k},-a\sqrt{q/k},aq/\sqrt{k},-aq/\sqrt{k},k/aq,kq^n,q^{-n};q)_j}{(q,\sqrt{a},-\sqrt{a},\sqrt{kq},-\sqrt{kq},\sqrt{k},-\sqrt{k},a^2q^2/k,aq^{-n+1}/k,aq^{n+1};q)_j}q^j\\&\quad - \sum\limits_{j = 0}^{\infty}\dfrac{(a,q\sqrt{a},-q\sqrt{a},a\sqrt{q/k},-a\sqrt{q/k},a/\sqrt{k},-aq/\sqrt{k},k/a,kq^n,q^{-n};q)_j}{(q,\sqrt{a},-\sqrt{a},\sqrt{kq},-\sqrt{kq},q\sqrt{k},-\sqrt{k},a^2q/k,aq^{-n+1}/k,aq^{n+1};q)_j}q^j
\\
&= \sum\limits_{j = 1}^{\infty}\dfrac{(a,q\sqrt{a},-q\sqrt{a},a\sqrt{q/k},-a\sqrt{q/k},-aq/\sqrt{k},kq^n,q^{-n};q)_jq^j}{(q,\sqrt{a},-\sqrt{a},\sqrt{kq},-\sqrt{kq},-\sqrt{k},aq^{-n+1}/k,aq^{n+1};q)_j}\\&\quad\times\left(\dfrac{(aq/\sqrt{k},k/aq;q)_j}{(\sqrt{k},a^2q^2/k;q)_j} - \dfrac{(a/\sqrt{k},k/a;q)_j}{(q\sqrt{k},a^2q/k;q)_j}\right)
\\
&= \dfrac{(1-a)(\sqrt{k}-a^2q/k+a/\sqrt{k}-k/aq)}{(1-\sqrt{k})(1-a^2q/k)}
\end{align*}
\newpage
\begin{align*}
\\&\quad\times\sum\limits_{j = 0}^{\infty}\dfrac{(aq,q\sqrt{a},-q\sqrt{a},a\sqrt{q/k},-a\sqrt{q/k},-aq/\sqrt{k},kq^n,q^{-n};q)_{j+1}}{(q;q)_j(\sqrt{a},-\sqrt{a},\sqrt{kq},-\sqrt{kq},q\sqrt{k},-\sqrt{k},a^2q^2/k,aq^{-n+1}/k;q)_{j+1}}\\&\qquad\qquad\times\dfrac{(aq/\sqrt{k},k/a;q)_j}{(aq^{n+1};q)_{j+1}}q^{j+1}
\\
&= \dfrac{q(1-aq)(1-aq^2)(1-kq^n)(1-q^{-n})(\sqrt{k}-a^2q/k+a/\sqrt{k}-k/aq)}{(1-k)(1-kq)(1-q\sqrt{k})(1-aq/\sqrt{k})(1-aq^{-n+1}/k)(1-aq^{n+1})}\\&\quad\times\sum\limits_{j = 0}^{\infty}\dfrac{(aq^2,q^2\sqrt{a},-q^2\sqrt{a},aq\sqrt{q/k},-aq\sqrt{q/k},aq/\sqrt{k},-aq^2/\sqrt{k},k/a;q)_j}{(q,q\sqrt{a},-q\sqrt{a},q\sqrt{kq},-q\sqrt{kq},q^2\sqrt{k},-q\sqrt{k},a^2q^3/k;q)_j}\\&\qquad\qquad\times\dfrac{(kq^{n+1},q^{-n+1};q)_j}{(aq^{-n+2}/k,aq^{n+2};q)_j}q^j
\\
&= \dfrac{q(1-aq)(1-aq^2)(1-kq^n)(1-q^{-n})(\sqrt{k}-a^2q/k+a/\sqrt{k}-k/aq)}{(1-k)(1-kq)(1-q\sqrt{k})(1-aq/\sqrt{k})(1-aq^{-n+1}/k)(1-aq^{n+1})}\\&\qquad\times B_{n-1}(aq^2,kq^2;q).
\end{align*}\\
Thus, we have\\
\begin{align}\label{eq217}
A_n(a,k;q) - B_n(a,k;q) = \psi_n(a,k;q)B_{n-1}(aq^2,kq^2;q),
\end{align}\\
where\\
\begin{align*}
&\quad \psi_n(a,k;q)\\ &= \dfrac{q(1-aq)(1-aq^2)(1-kq^n)(1-q^{-n})(\sqrt{k}-a^2q/k+a/\sqrt{k}-k/aq)}{(1-k)(1-kq)(1-q\sqrt{k})(1-aq/\sqrt{k})(1-aq^{-n+1}/k)(1-aq^{n+1})}.    
\end{align*}\\
Similarly, after tedious computations, it can also be shown that\\
\begin{align}\label{eq218}
B_n(a,k;q) - A_{n-1}(a,k;q) = \chi_n(a,k;q)B_{n-1}(aq^2,kq^2;q),    
\end{align}\\
where\\
\begin{align*}
\chi_n(a,k;q) &= \dfrac{(1-aq)(1-aq^2)(1-kq^n)}{(1-k)(1-q\sqrt{k})(1-kq)(1-aq^n)(1-aq^{n+1})}\\&\qquad\times\dfrac{\Tilde{\chi}_n(a,k;q)}{(1-k^2q^{n-2}/a^2)(1-kq^{n-1}/a)},    
\end{align*}\\
where\\
\begin{align*}
&\quad \Tilde{\chi}_n(a,k;q)\\ &= -k^2q^{n-2}/a^2 - k^2q^{n-1}/a^2 + k^4q^{2n-3}/a^4 - aq^n + k^2q^{2n-1}/a - k^4q^{3n-3}/a^3\\&\quad - \sqrt{k} + k^2q^{n-2}\sqrt{k}/a^2 - k^4q^{2n-3}\sqrt{k}/a^4 + aq^n\sqrt{k} - k^2q^{2n-2}\sqrt{k}/a\\&\quad - k^2q^{2n-1}\sqrt{k}/a + kq^{n-1}/a + k^2/a^2q - k^3q^{n-2}/a^3 + kq^{n-1} - k^3q^{n-2}/a^2\\&\quad + k^4q^{2n-3}/a^3 + q^n\sqrt{k} - kq^{2n-1}\sqrt{k}/a + k^3q^{2n-2}\sqrt{k}/a^3 - kq^{2n-1}\sqrt{k}\\&\quad + k^2q^{3n-2}\sqrt{k}/a + k^3q^{2n-2}\sqrt{k}/a^2.
\end{align*}\\ Now, let $A_n^{\prime}(a,k;q)$ and $B_n^{\prime}(a,k;q)$ denote the right-hand sides of \eqref{eq113} and \eqref{eq114} respectively. Then, we have\\
\begin{align*}
&\quad A_n^{\prime}(a,k;q) - B_n^{\prime}(a,k;q)\\ &= \dfrac{(aq,k^2/a^2q;q)_n}{(k,k/a;q)_n} - \dfrac{(aq,\sqrt{k},k^2/a^2;q)_n}{(k,k/a,q\sqrt{k};q)_n}
\\
&= \dfrac{(aq,k^2/a^2q;q)_n}{(k,k/a;q)_n} - \dfrac{(1-\sqrt{k})(aq,k^2/a^2;q)_n}{(1-q^n\sqrt{k})(k,k/a;q)_n}
\\
&= \dfrac{(1-q^n)(\sqrt{k}-k^2/a^2q)(aq;q)_n(k^2/a^2;q)_{n-1}}{(1-q^n\sqrt{k})(k,k/a;q)_n}
\\
&= \dfrac{(1-q^n)(\sqrt{k}-k^2/a^2q)(kq^2,k/a;q)_{n-1}(aq;q)_n(aq^3,q\sqrt{k},k^2/a^2;q)_{n-1}}{(1-\sqrt{k})(k,k/a;q)_n(aq^3;q)_{n-1}(kq^2,k/a,q^2\sqrt{k};q)_{n-1}}
\\
&= \dfrac{q(1-aq)(1-aq^2)(1-kq^n)(1-q^{-n})(\sqrt{k}-a^2q/k+a/\sqrt{k}-k/aq)}{(1-k)(1-kq)(1-q\sqrt{k})(1-aq/\sqrt{k})(1-aq^{-n+1}/k)(1-aq^{n+1})}\\&\qquad\times B_{n-1}^{\prime}(aq^2,kq^2;q).
\end{align*}\\
Thus, we have\\
\begin{align}\label{eq219}
A_n^{\prime}(a,k;q) - B_n^{\prime}(a,k;q) = \psi_n^{\prime}(a,k;q)B_{n-1}^{\prime}(aq^2,kq^2;q),
\end{align}\\
where\\
\begin{align*}
\psi_n^{\prime}(a,k;q) &= \dfrac{q(1-aq)(1-aq^2)(1-kq^n)(1-q^{-n})}{(1-k)(1-kq)(1-q\sqrt{k})(1-aq/\sqrt{k})}\\&\qquad\times\dfrac{(\sqrt{k}-a^2q/k+a/\sqrt{k}-k/aq)}{(1-aq^{-n+1}/k)(1-aq^{n+1})}.    
\end{align*}\\
Again, it can be shown that\\
\begin{align}\label{eq220}
B_n^{\prime}(a,k;q) - A_{n-1}^{\prime}(a,k;q) = \chi_n^{\prime}(a,k;q)B_{n-1}^{\prime}(aq^2,kq^2;q),
\end{align}\\
where\\
\begin{align*}
\chi_n^{\prime}(a,k;q) &= \dfrac{(1-aq)(1-aq^2)(1-kq^n)}{(1-k)(1-q\sqrt{k})(1-kq)(1-aq^n)(1-aq^{n+1})}\\&\qquad\times\dfrac{\Tilde{\chi}_n^{\prime}(a,k;q)}{(1-k^2q^{n-2}/a^2)(1-kq^{n-1}/a)},    
\end{align*}\\
where\\
\begin{align*}
&\quad \Tilde{\chi}_n^{\prime}(a,k;q)\\ &= -k^2q^{n-2}/a^2 - k^2q^{n-1}/a^2 + k^4q^{2n-3}/a^4 - aq^n + k^2q^{2n-1}/a - k^4q^{3n-3}/a^3\\&\quad  - \sqrt{k} + k^2q^{n-2}\sqrt{k}/a^2 - k^4q^{2n-3}\sqrt{k}/a^4 + aq^n\sqrt{k} - k^2q^{2n-2}\sqrt{k}/a\\&\quad  - k^2q^{2n-1}\sqrt{k}/a + kq^{n-1}/a + k^2/a^2q - k^3q^{n-2}/a^3 + kq^{n-1} - k^3q^{n-2}/a^2\\&\quad + k^4q^{2n-3}/a^3 + q^n\sqrt{k} - kq^{2n-1}\sqrt{k}/a + k^3q^{2n-2}\sqrt{k}/a^3 - kq^{2n-1}\sqrt{k}\\&\quad + k^2q^{3n-2}\sqrt{k}/a + k^3q^{2n-2}\sqrt{k}/a^2.    
\end{align*}\\
Now, observe that \eqref{eq217} and \eqref{eq218} together with the initial conditions $A_0(a,k;q)$ $=$ $B_0(a,k;q) = 1$ uniquely define both $A_n(a,k;q)$ and $B_n(a,k;q)$. Similarly, \eqref{eq219} and \eqref{eq220} together with the initial conditions $A_0^{\prime}(a,k;q)$ $=$ $B_0^{\prime}(a,k;q) = 1$ uniquely define both $A_n^{\prime}(a,k;q)$ and $B_n^{\prime}(a,k;q)$. Also, note that the recurrences \eqref{eq217} and \eqref{eq219} are the same and the recurrences \eqref{eq218} and \eqref{eq220} are the same. Hence,\\
\begin{align*}
A_n(a,k;q) &= A_n^{\prime}(a,k;q)
\end{align*}
and
\begin{align*}
B_n(a,k;q) &= B_n^{\prime}(a,k;q).
\end{align*}\\
Thus, we have a proof of the two ${}_{10}W_9$ identities due to Andrews-Berkovich (Theorem \ref{thm19}) using Pfaff's method.\qed\\

\section*{Acknowledgments}
The author would like to thank George E. Andrews for advising him to look at Pfaff's method as a comment after his talk titled ``On ${}_5\psi_5$ identities of Bailey'' at the Number Theory Seminar in the Department of Mathematics, University of Florida on March 21, 2023. The author would also like to thank his Ph.D. advisor, Alexander Berkovich, for suggesting him to prove Theorem \ref{thm19} and for his continuous support and encouragement. The author would also like to thank the anonymous referee for valuable comments and suggestions.

\nocite{*}
\bibliographystyle{amsplain}


\end{document}